\documentclass[10pt]{article}
\usepackage{stix2} 
 \usepackage[cal=boondoxo]{mathalfa} 
   \usepackage{amsmath,amsthm,amscd,graphicx}
 \usepackage[pagebackref]{hyperref} 
 \usepackage{color}

\newtheorem{thm}{Theorem}[section]
\newtheorem*{thm*}{Theorem}
\newtheorem*{corr*}{Corollary}
\newtheorem{lemma}[thm]{Lemma}

\newtheorem*{prop*}{Proposition}
\newtheorem{corr}[thm]{Corollary}
\theoremstyle{definition}

\newtheorem*{dfn*}{Definition}
\newtheorem{exmple}[thm]{Example}
\newtheorem*{exmple*}{Example}
\newtheorem{exmples}[thm]{Examples}

\newtheorem*{conj*}{Conjecture}
\theoremstyle{remark}

\newtheorem{rmk}[thm]{\textit{Remark}}
\newcommand{\bbF}{{\mathbb{F}}}

\newcommand{\bP}{{\mathbb{P}}}
\newcommand{\bQ}{{\mathbb{Q}}}
\newcommand{\bR}{{\mathbb{R}}}
\newcommand{\bZ}{{\mathbb{Z}}}

 \newcommand{\comp}{\raise1pt\hbox{{$\scriptscriptstyle\circ$}}}
 \def\disc#1{\operatorname{\rm disc}(#1)}%
 \def\lset{\{}  
\def\rset{\}}  
\def\set#1{\lset#1\rset} 

\def\hom{\text{\rm Hom}}

\def\mapright#1{\mathop{\vbox{\ialign{
                                ##\crcr
    ${\scriptstyle\hfil\;\;#1\;\;\hfil}$\crcr
 \noalign{\kern2pt\nointerlineskip}
    \rightarrowfill\crcr}}\;}}

\def\half{\frac 12}

 \def\disc#1{\text{\rm disc}\left(#1\right)}%

\def\la{\langle}
\def\ra{\rangle}
\def\rank#1{\operatorname{rank}({#1})}
\begin{document}

\title{A note on the primitive cohomology  lattice  of a projective surface}
\author{Chris Peters\\
Technical University Eindhoven}
\date{Revised January 2022}

\maketitle 
   \abstract{The isometry class of the intersection form of a compact complex surface 
   can be easily determined from complex-analytic invariants. For projective surfaces the primitive lattice
   is another naturally occurring lattice. The goal of this note is to show that it can be determined from
   the intersection lattice and the self-intersection of a  primitive ample class, at least when the
   primitive lattice is indefinite.  Examples include the Godeaux surfaces, the Kunev surface and a specific Horikawa surface.  There are also some results concerning (negative) definite primitive lattices, especially for canonically polarized surfaces of general type. 
   \\
   \textbf{Keywords:} Complex projective surfaces, primitive intersection lattice\\
 \textbf{MSC:} 14J80, 32J15, 57N65
     } 
   
\section{Introduction}

The intersection form of a compact connected orientable $4n$-dimensional  manifold $X$ is the   bilinear, symmetric form on $\mathsf H_X= H^{2n}(X,\bZ)/\text{(torsion)}$
given by cup product. By Poincar\'e duality this form is  \textbf{\emph{unimodular}}, that is, its Gram matrix has determinant $\pm1$. 
The pair  consisting of $\mathsf H_X$ and the intersection pairing is called the \textbf{\emph{intersection lattice}}  of $X$.

If $X\subset \bP^N $ is a smooth compact complex  manifold with hyperplane section $H$,  the orthogonal complement of  the class of $H^n$ in 
   $ H^{2n}(X,\bZ)$   
is called the (middle) \textbf{\emph{primitive cohomology}}, 
denoted $\mathsf P_X $.  Precise knowledge of this lattice and its group of isometries  
turns out to be useful, especially for arithmetic questions.
This motivates interest in the  main result of this note which deals with the case of surfaces (=Theorem~\ref{thm:main}):   

\begin{thm*}  Let $X$ be a complex projective surface  with $p_g(X)\not=0$ and let $c\in \mathsf H_X$ be  a primitive representative of an ample divisor.
 Then     the isometry class of the lattice $c^\perp$ is uniquely determined by the following data:
 \begin{enumerate}
\item the  triple  $( b_1(X), c_1^2(X), c_2(X))$ of topological invariants,
\item whether or not  $c$ is characteristic
\item the self-intersection of $c$.
\end{enumerate}
  \end{thm*}

The above result implies in particular that for a given surface $X$   the primitive lattice  does not not depend on 
the particular choice of the projective embedding of $X$, but only on the  degree of $X$.
The proof  of the theorem uses  firstly  Nikulin's reformulation of the classical classification results on integral quadratic forms 
in terms of   the discriminant quadratic form and, secondly,  on a fine analysis of the type of intersection lattices occurring  for projective surfaces based on the Enriques classification.
This result  is effective as   illustrated
for surfaces with small $c_1^2$, e.g.\   for some  Horikawa surfaces.   See Examples~\ref{ex:calcs}.
 
  The assumption $p_g(X)\not=0$ is equivalent to $\mathsf P_X$ being indefinite, a prerequisite for applying Nikulin's results. However in the definite situation
 one can  in several instances still determine  the isometry class of the primitive intersection lattice making use of a series of investigations by G. Watson~\cite{wats,wats1, wats2,wats3,wats4,wats5,wats6,wats7}. 
 See Remark~\ref{rmk:OnIndefLats}.
 
 \begin{rmk} Primitive cohomology   plays a central role in Hodge theory since the Hodge decomposition together with the  intersection pairing gives  
$\mathsf P_X$ the  structure of a polarized pure Hodge structure of weight $2n$.  
To explain why this is the case,  consider an embedding  $X\subset \bP^N$.
The Hodge structure on the middle primitive cohomology in smooth families 
$\set{X_s}_{s\in S}$ of smooth varieties embedded in the same $\bP^N$  gives rise to a period map $S\to \Gamma\backslash D$ 
where $D$ is a suitable period domain and where $\Gamma$, the (maximal) monodromy group,  is
the isometry group of  the primitive lattice of a fibre $X_s$ (all such groups are isomorphic). More precisely, since monodromy preserves the polarization,
$\Gamma$ is the subgroup of the isometry group of $P=\mathsf P_{X_s}$ inducing the identity on the discriminant group $P^*/P$.
 \end{rmk}

\subsection*{Acknowledgement} I want to thank G.\ Pearlstein for asking me how to compute the primitive cohomology
of the Horikawa surface with $c_1^2=1, c_2=35$ which inspired me to write this note.   Thanks also go to 
M.\  Sch\"utt for a careful reading of a first draft.

 \subsection*{Conventions and Notations} \begin{itemize}
\item A lattice is a free $\bZ$-module of finite rank equipped with a non-degenerate symmetric bilinear integral form which is denoted 
 with a dot. 
 \item A rank one lattice $\bZ e$ with $e.e= a$ is denoted $\la a\ra$, orthogonal direct sums by $\operp$.
  Other standard lattices are the hyperbolic plane $U$, and the root-lattices $A_n,B_n$ ($n\ge 1$), $D_n$, $n\ge 4$. and $ E_n$, $n=6,7,8$.
  Their $p$-adic localizations will be denoted by the same symbol. 
  More details are given below in Section~\ref{sec:lats}. 
  
  \item    If one replaces the form on the lattice $L$ by $m$-times the form, $m\in \bZ$, 
  this  scaled lattice is denoted   $L(m)$. 
  
  \item An inner product space over a field $k$ is a $k$-vector space equipped with a non-degenerate symmetric bilinear form over $k$. It will likewise be denoted with a dot. 
 \item The signature of a    non-degenerate symmetric bilinear integral form  $b$  is denoted by $(b^+ ,b^- )$ 
 and the index by $\tau=b^+-b^- $.
The signature of  the intersection lattice $\mathsf H_X= H^{2n}(X,\bZ)/\text{torsion}$, $X$ a compact connected orientable $4n$-dimensional  manifold,   will be denoted by $\tau(X)$.     If $X$ is projective, its "primitive cohomology" is  the integral primitive cohomology (classes of $\mathsf H_X$ orthogonal to an ample class) and is denoted by $\mathsf P_X$.

\end{itemize}

\section{On lattices} \label{sec:lats}

 \subsection*{Unimodular lattices} As is well known (cf. \cite{milnor58,milnor})  if a unimodular form is indefinite, its isometry class is uniquely determined by the signature
  and type of the form.  The \emph{type} of  a bilinear symmetric  form  by definition is  \emph{even} or \emph{odd}. 
 Being  even means that $x.x$ is even for all elements $x$ of the lattice
  and odd otherwise.  The results from loc.\ cit. state that odd unimodular forms are diagonalizable over the integers.
  This is evidently not the case for unimodular even forms. Instead these are orthogonal sums of
   three  building blocks, the hyperbolic plane $U$, the positive definite root lattice $E_8$,
   and its negative $E_8(-1)$. The first has rank two and has a basis $\set{e,f}$
  for which $e.e=f.f=0$ and $e.f=1$. The root lattice $E_8$ has rank $8$ with form given in the basis  by the Coxeter matrix for the root lattice $E_8$, that is by 
  \[\begin{pmatrix}
  2 & 0 & -1 & 0 & 0 & 0 & 0 & 0   \\
  0 & 2 & 0 & -1 & 0 & 0 & 0 & 0   \\
  -1 & 0 & 2 & -1 & 0 & 0 & 0 & 0   \\
   0 & -1 & -1 & 2 & -1 & 0 & 0 & 0   \\
   0 & 0 &0 & -1 & 2 & -1 & 0 & 0   \\
   0 & 0 &0 & 0& -1 & 2 & -1 & 0   \\
   0& 0 & 0 &0 & 0& -1 & 2 & -1 \\
   0& 0& 0 & 0 &0 & 0& -1 & 2
  \end{pmatrix} .
  \]
 It turns out that every indefinite even unimodular form is isometric to $\operp^s U\operp \operp^t E_8$,
 a lattice of  index $t\ge 0$,  or to $\operp^s U\operp \operp^t E_8(-1)$ if the index equals   $-t<0$.

 For definite forms the situation is more complicated. The number of non-isometric lattices grows rapidly with the rank. See e.g. \cite[Ch.\ IV \S~2.3]{serre}.

\subsection*{Characteristic elements}  To test whether the form on a lattice $L$ is even or odd, one makes use of a \textbf{\emph{characteristic element}} $c\in L$. By definition it has the property that $c.x+x.x$ is even for all $x\in L$. Such characteristic elements exist if the discriminant of $L$ is odd as one easily sees by reduction modulo $2$.
 In fact, characteristic classes exist for inner product spaces over the field $\bbF_2$.   Of course, if $c\in L$ is not isotropic and $L$ is even, then  $c^\perp$ is an even lattice, but this holds   also if $c$
 is characteristic in an odd lattice $L$. For later use  I  set this apart:
 
 \begin{lemma}\label{lem:OrtChar} If $L$ is a lattice with odd discriminant and  $c\in L$ not isotropic, i.e.\  $c\cdot c\not=0$, then $c^\perp$ is an even  lattice if and only if $c$ is a characteristic element.
 \end{lemma}
 \begin{rmk} An odd unimodular indefinite lattice being diagonalizable, the reader may be surprised that  it  can have   unimodular even sublattices.
 That this is indeed the case can be illustrated    with the lattice  $L=\la 1\ra \operp \la 1\ra \operp \la -1\ra$. The
 basic observation is that $L$ is isometric to $\la 1\ra \operp U$.
 Explicitly, if $\set{e_1,e_2,e_3}$ is an orthogonal   basis for $L$, then $c=2c'$, $c'=e_1+e_2+e_3$ is a characteristic element
 with $c'.c'=1$ and $c^\perp$ is the lattice with basis $\set{e_1+e_3,e_2+e_3}$ isometric to $U$. 
 \end{rmk}

 \subsection*{Discriminant forms and the genus}
 
 Let $L$ be a lattice. We  recall the concept of discriminant group and discriminant form.
  Remark that the pairing on $L$ extends to a $\bQ$-bilinear 
 pairing  on $L\otimes\bQ$ and induces   the $\bQ/\bZ$-valued form  on  the  discriminant group   $A(L)=L^*/L$,   $L^*=\hom_\bZ(L,\bZ)$ given by
 \[ 
b_L  :   A(L)\,  \times  \, A(L)  \to \bQ/\bZ,   \quad  \bar x . \bar y  \mapsto  x. y  \bmod \bZ \,\text{(\textbf{\emph{discriminant bilinear form}})} .
 \] 
 Even lattices come with an integral quadratic form $q$ given by $q(x)= \half b(x,x)$ and for these  one considers   a finer invariant, the 
 \textbf{\emph{discriminant quadratic form}}
 \[
 q_L : A(L) \to  \bQ/\bZ,   \quad \bar x \mapsto   {q(x)} \bmod \bZ.
 \]
 The  discriminant form is completely local in the sense that it decomposes into $p$-primary forms where $p$ is a prime dividing the discriminant.
 More precisely, it is the orthogonal sum of the discriminant forms of the localizations $L_p=L\otimes\bQ_p$ and so it ties in with the \textbf{\emph{genus}} of the lattice,  i.e.\ the  set of isometry classes  $\set{L_p}_{p\text { prime}}$ together with $L\otimes\bR$. A celebrated result  of V. Nikulin~\cite[Cor. 1.16.3]{nikulin} emphasizes the role of the discriminant form in determining  the genus:
 
 \begin{thm*} The genus of  non-degenerate lattice is completely determined by its type, rank,   index and the discriminant form.
 \end{thm*}
 
It is well known that the number of isometry classes in a genus is finite. It is also called the \textbf{\emph{class number}} of the genus.

 For applications in geometry it is important to have a criterion for class number $1$ lattices. This is  often the case
 in the indefinite situation as stated by another result due to V.\ Nikulin \cite[1.13.3 and 1.16.10]{nikulin} and M.\ Kneser \cite{kneser}:

  \begin{thm}
 Let $L$ be a  non-degenerate \textbf{indefinite} lattice of rank $r$. Its class number is $1$ in the following instances:
 \begin{enumerate}
\item  In case $L$   is even and  the discriminant group  of $L$
can be generated by    $\le r-2$ elements.  Hence, in this case
$L$  is uniquely determined by  its  rank,   index   and the discriminant  quadratic   form.
\item In case $L$  is odd, and the discriminant group  of $L$
can be generated by     $\le r-3$ elements. 
Hence, in this case
$L$  is uniquely determined by its  rank,   index   and the discriminant
bilinear form.
\end{enumerate}

\label{thm:kneser}
\end{thm}

These results   will be  in particular applied to  primitive 
   sublattices  of $L$, i.e.\   sublattices  $S$ such that $L/S$ is free of torsion. 
In case $S$ is well understood, one can say much about its orthogonal complement:

\begin{lemma} \label{lem:SandT} Let $S$ be a primitive non-degenerate sublattice of $L$ and $T=S^\perp$ then  $\disc S=\pm \disc {S^\perp}$  and $(A(S),b_S)$ is isometric to $(A(T),-b_T)$.
\end{lemma}
For proofs, see e.g.\  \cite{K3lects,kondoK3}.
 
\subsection*{Intersection lattices}

Lemma~\ref{lem:SandT} has the following implication for intersection lattices:

\begin{corr} 
\label{cor:Main} Let $X$ be a compact connected orientable $4n$-dimensional  manifold $X$ with indefinite intersection form and let $c\in \mathsf H_X$ be primitive  with $c.c\not=0$.
If $\mathsf H_X$ is even assume  that $b_n(X)\ge 4$    and if $\mathsf H_X$ is odd and $c$ is not characteristic, assume that  $b_n (X) \ge 5$.
Then the isometry class of $c^\perp$ is uniquely determined by the signature $(b^+,b^-)$ of $\mathsf H_X$ and the  integer   $ c.c $.

\end{corr}
\begin{proof} The discriminant form of $\bZ.c$ equals $\la 1/(c.c)\ra$ and by Lemma~\ref{lem:SandT} the discriminant form for $T:=c^\perp$
equals $ - \la 1/(c.c)\ra$ and, in particular, is a torsion form on a length  one group. The assumptions imply that $1\le \rank T- 2$ in the even case, and
$1\le \rank T- 3$ in the odd case. Since $T$ is odd if and only if $c$ is not characteristic, the statement follows.
\end{proof}

Assume now that $X$ is a compact  orientable $4$-dimensional  manifold with intersection lattice $\mathsf H_X$.  The second Stiefel--Whitney class $w_2$ is a  characteristic class for  the inner product space $H^2(X,\bbF_2)$.  To pass to integral cohomology one uses the
reduction mod $2$  map, induced by the natural projection $\bZ\to \bZ/2\bZ$:
\begin{align} \label{eqn:OnMod2Red}
 \rho_2 :H^2(X,\bZ) \to  H^2(X,\bZ/2\bZ)  .
 \end{align}
 Any lift of $w_2$ under $\rho_2$ is an  integral  characteristic element since the
intersection pairing is compatible with reduction modulo $2$.
In the special case where $X$ is a compact almost complex manifold of complex dimension $2$, there is a canonical choice for a lift, namely
the first Chern class $c_1$. We note a simple consequence:

\begin{lemma} The intersection pairing on a compact almost-complex surface $X$   is even if   $c_1(X)$ is divisible by $2$ in integral cohomology.
 The converse is true if $H_1(X,\bZ)$ is free of $2$-torsion.
 
If $c^2_1(X)\not=0$, and $c_1 (X)= kc$ with $c$ primitive, then 
 the lattice  $c_1(X)^\perp\subset H^2(X,\bZ)$ is a non-degenerate  even lattice of discriminant $\pm c.c$.
\end{lemma}

\begin{proof} The preceding remarks show that if $c_1(X)$ is divisible by $2$ in cohomology, $x .  x $ is even for all $x\in H^2(X,\bZ)$.
For the converse, consider the long exact sequence associated to $
0 \to \bZ \mapright{ \times 2} \bZ \to \bZ/2\bZ\to 0
$ and use that the intersection pairing on $H^2(X,\bbF_2)$ is non-degenerate.   Here 
  surjectivity of the map $\rho_2$  (cf. \eqref{eqn:OnMod2Red}) is used which follows  since  by Poincar\'e-duality, $H^3(X,\bZ)\simeq H_1(X,\bZ)$ -- which  has no $2$-torsion by assumption.

The penultimate  assertion is also clear since $c_1. x+ x.x=x.x$ is even for all $x\in c_1(X)^\perp$. The assertion about the discriminant is a special case of Lemma~\ref{lem:SandT}.
\end{proof}

 \section{Complex algebraic  surfaces} \label{sec:surfs}
 
 \subsection*{Invariants}
 
 Let $X$ be a compact complex projective surface. It is well known that the Chern numbers, a priori complex invariants,  are in fact  (oriented) topological
invariants. This is clear for $c_2(X)$  since it can be identified with the Euler number $e(X)$.
To see that $c_1^2(X)$ is a topological invariant, one  invokes  a deep theorem, the index theorem (\cite[Thm. 8.2.2]{hirz}):

\begin{thm}[Index theorem -- special case]   \index{index theorem}
\label{thm:Index}
For a compact differentiable $4$-manifold  $X$ admitting a complex structure, the index  $\tau(X)$  satisfies
\[
\tau(X) = \frac 13 (c_1^2(X)   - 2 c_2(X)) .
\] 
\end{thm}

For algebraic surfaces the Hodge decomposition gives two more invariants for $X$, namely $q(X)= \half b_1(X)$ and $p_g(X)= \dim H^{2,0}(X)$. 
In particular, $q $ is a topological invariant. Because of  \textbf{\emph{Noether's formula}} \cite[p. 26]{4authors},   
\begin{align} \label{eqn:Noether}
 \chi(X):=1-q(X) +p_g (X) 
= \frac{1}{12}(c_1^2(X)+c_2(X)),
\end{align}
also $p_g$  is a topological invariant. 

Recalling that since $c_1$ is a characteristic element for the intersection lattice, these observations  make it possible to
determine the isometry class of $\mathsf H_X$ from the type of $c_1$ together with the integer invariants $c_1^2$ and $c_2$.

\begin{exmple} A \textbf{\emph{K3 surface}} by definition is a surface with $b_1=0$  and trivial canonical bundle and so $c_1=0$ and $p_g=1$, $q=0$ implying 
  $2= \frac{1}{12} c_2$.  Hence $b_2=24-2=22$. The index theorem gives $\tau= \frac 13 (-48)= -16$ and since the intersection lattice is even,
it is isometric to $\operp^3 U \operp \operp^2 E_8(-1)$.

An \textbf{\emph{Enriques surface}} has $p_g=q=0$ while  $c_1$ is $2$-torsion. A similar reasoning shows that $U\operp E_8(-1)$ is its intersection lattice.

\end{exmple}

For algebraic surfaces (and more generally for compact K\"ahler surfaces) there is a characterization of the signature in terms of Hodge numbers:

\begin{lemma}[\protect{\cite[Thm. IV.2.6]{4authors}}] Let $X$ be a compact K\"ahler surface. Then the  signature  of $X$ equals $(2p_g(X)+1,h^{1,1}(X) -1)$ where $h^{1,1}(X)=\dim H^{1,1}(X)$. \label{lem:OnSign} 
\end{lemma}

 \subsection*{Surface classification}

 I also make use of the Enriques classification of surfaces. The notion of a minimal surface plays an essential role. All surfaces are obtained from these by
 repeated blowing up in points. In the present context it is important to recall how the intersection lattice changes under a blow-up.
 Since blowing up $X$ in a point does not affect $H^i$, $i\not=2$ and replaces $H^2(X)$ by $H^2(X)\oplus \bZ$, where the summand $\bZ$ is generated by the exceptional
 curve which has self-intersection $-1$, one has:
 
 \begin{lemma} \label{lem:OnBlop} Let $X$ be a compact complex surface and let $\widetilde X  $ be the surface obtained by blowing up  $X$ in a point.
 Then $\mathsf H_{\widetilde X}= \mathsf H_X\operp \la -1 \ra$. In particular, the intersection lattice of a non-minimal surface is odd.

 Moreover $c_1^2(\widetilde X)= c_1^2(X)-1$, $c_2(\widetilde X)=c_2(X)+1$ and $\tau(\widetilde X)=\tau(X)-1$.
 \end{lemma}
 
In the Enriques classification -- besides the already mentioned classes (K3 surfaces, Enriques surfaces) -- some other
classes appear.  Firstly the rational and ruled surfaces which by definition are obtained from the projective plane, respectively a minimal ruled surface by repeatedly blowing up and blowing down. 
Then there are the elliptic surfaces  which by definition  admit  a holomorphic map  onto  a curve such that the  general fibre is an
elliptic curve.  Among these are some ruled surfaces, the Enriques surfaces and some K3 surfaces. Next,    there are  so-called bi-elliptic or hyperelliptic surfaces and, finally,  the large class of properly elliptic surfaces which by definition have Kodaira dimension $1$. 
The surfaces with Kodaira dimension $2$ are called "surfaces of general type". Together these exhaust the classification (see e.g.\ \cite{CAS}).  Summarizing, replete with invariants, one has:
 
 \begin{thm}[Enriques  classification]\label{thm:EnrClass}
Every minimal   complex projective   surface belongs to exactly one of the
following classes ordered according   to their Kodaira dimension $\kappa$:
 
\begin{center}
\begin{tabular}{|c | l | c| c |c| c|c|}
 \hline
$\kappa$ & Class  								&   $b_1$ & $p_g$ & 	&$c_1^2$    	& $c_2$ \\
 \hline
$ -\infty$  &  minimal rational  surfaces 			       & $ 0$   & $0$ & 		& $8$ or  $9$  & $4$ or $3$   \\
		 &  ruled surfaces of genus $>0$   			&  $ 2g$ 	 & $0$ &  & $8(1-g)$&    $4(1-g)$ \\
		\hline
$0$ &	Two-dimensional tori 						&     $4$  & $2$ &  & $0$  & $0$ \\
	         &	K3 surfaces       			 	        & $0$  & $1$ &  &   $0$ &  $    24$  \\
	         & Enriques surfaces     					 & $0$  & $0$ &  & $0$   & $12$ \\
	         & bielliptic surfaces 					      &	 $2$ & $0$ &  & $0$  & $0$ \\
	        \hline
$ 1$	  &    minimal properly elliptic surfaces   		   &  		 &   & 	  & $0$ & $\ge 0$  \\
 \hline
  $  2$  &    minimal surfaces of general type  					&   &   &   & $>0$ & $>0$\\
   \hline
\end{tabular}
\end{center} 
 \end{thm}

In the next section one considers indefinite primitive lattices. 
Here I discuss the   -- rather small -- list of   surfaces having definite primitive  lattices. First of all, these cannot be positive definite:
\begin{lemma} Let $X$ be a  complex projective surface. Then $H^2(X,\bR)$ is positive definite
if and only if $b_2(X)=1$ and so  $\mathsf P_X\not=0$ cannot be  positive definite. \label{lem:ExeptMain1}
\end{lemma}
\begin{proof} Lemma~\ref{lem:OnSign} implies that the signature of a primitive lattice is $(2p_g,  h^{1,1}-1)$ which is
positive definite precisely if  $b_2=\tau= 2p_g+1$. Moreover, $X$ is minimal of general type and one finds $c_1^2=10p_g-8q+9$ and $c_2=2p_g-4q+3$.
Now invoke  the Bogomolov--Miyaoka--Yau inequality   (cf.\ \cite[\S VII.4]{4authors}), stating
\begin{align} \label{eqn:BMY}
c_1^2-3c_2\le 0,
\end{align}
 which gives $ 4p_g+4q \le 0$ and so  $p_g=q=0$.
But then $b_2=1$ which forces $\mathsf P_X=0$.  
\end{proof}

%
 
Secondly, as to negative definite $\mathsf P_X$, by Lemma~\ref{lem:OnBlop} one may restrict to minimal surfaces
and hence, inspecting the table from Theorem~\ref{thm:EnrClass}, one sees:

\begin{lemma} \label{lem:ExeptMain2}
 Let $X$ be complex projective surface  with $\mathsf P_X \not=0$   and negative definite. 
 Then $X$ is either   rational or ruled, a (possibly blown-up)  Enriques surface, an elliptic surface with $p_g=0$ or a surface of general type with $p_g=0$.
  \end{lemma}
 
 \begin{rmk} In the    definite  situation there might be more  isometry classes in the genus. 
 There are however instances where  the class number is exactly one. For minimal surfaces that are canonically polarized  and with $p_g=0$ this can be used to determine the primitive cohomology.
 See   the table in Remark~\ref{rmk:OnIndefLats}.
  \end{rmk}

In the next section one also needs the following result:

\begin{lemma} Let $X$ be a  complex projective surface  
with $b_2(X)\le 4$ and $p_g(X)=1$. Then $X$ is a minimal  algebraic  surface satisfying 
$c_1^2(X)= 3c_2(X)=18$, $q(X)=0$ (and so $b_2(X)=4$).\label{lem:Excepts} 
\end{lemma}

\begin{proof} Assume that $X$ is minimal elliptic. Since $p_g=1$ the surface is either K3 or properly elliptic. However, since $b_2\le 4$, the surface cannot be K3. So it is properly elliptic with invariants
 $c_1^2(X)=0$ and $c_2(X)= 12 (p_g(X)-q(X)+1) = 12(2-q)$.
On the other hand $c_2(X)=2-4q(X)+b_2(X)$ and so $4\ge b_2(X)\ge  2 p_g+1=3$ must be even and  hence 
 $b_2(X)=4$, but then $c_2(X)= 12(2-q(X))= 6-4q(X)$ which is impossible. If $X$ is not minimal, for  its minimal model we have  $b_2 \le 3$ and so it also does not exist

If $X$ is  of general  type, then from $c_2(X)= 2-4q(X)+b_2(X)>0$ one finds $q(X)=0$. Since $b_2(X)=3 ,4 $, from 
 $24 = c_1^2(X) + c_2(X)$   one finds that either $(c_1(X), c_2(X)) =(19 ,5)$ or $=(18, 6)$.
The inequality  \eqref{eqn:BMY} excludes the first possibility and then $b_2(X)=4$. If $X$ were not minimal and $\widetilde X$ its minimal model,
then $b_2(\widetilde X)=3$ which is excluded by the previous calculation.
\end{proof}

\begin{rmk} Since $X$ satisfies $c_1^2(X)= 3c_2(X)$,  by S.T.\ Yau's results \cite{yau}, its universal cover is the unit ball.
The existence of a surface with $p_g(X)=1, q(X)=0$ and $c_1^2(X)=18$ is not known.
These are of course far from simply connected. For simply connected surfaces the maximum $c_1^2$ seems to
be $12$ (G.\ Urzua, unpublished).
\end{rmk}

 \section{On primitive  intersection lattices of surfaces}
 
 \label{sec:primlaats}

 The main result is as follows.
 
 \begin{thm} Let $X$ be complex projective surface  whose primitive lattice $\mathsf P_X$  is  indefinite. Let $h\in \mathsf H_X$ be a primitive ample  class. Then
 \begin{enumerate}
\item  If \ $\mathsf H_X$ is   even, the
 isometry class of $\mathsf P_X$  is uniquely determined by  the triple  $( b_1(X), c_1^2(X), c_2(X))$ of topological invariants together with
 $h.h$.
 \item In  case $\mathsf H_X$ is odd, this depends in addition to   $h$ being characteristic or not: In case $h$ is characteristic,  $\mathsf P_X$ is even and otherwise it is odd.
If the latter occurs, one  assumes  in addition  that $b_2(X)\not= 4$.
\end{enumerate}
\label{thm:main}
 \end{thm}

 \begin{proof} This is a direct consequence of Corollary~\ref{cor:Main}.
Indeed,  since $\tau(X)=\frac 13 (c_1^2(X)-c_2(X))$, the index of $\mathsf P_X$ equals $\tau(X)-1$ and
 $\rank{\mathsf P_X}=b_2(X)-1=c_2(X)-2b_1(X)-1$. The result follows from  Corollary~\ref{cor:Main}  and Lemma~\ref{lem:Excepts}.
 Indeed, the latter result implies that $b_2(X)\ge 4$.  \end{proof}

  \begin{exmples}  \label{ex:calcs} 
 \begin{enumerate}
 
\item  For a   complex projective surface $X$  with $c_1^2(X)=1$ and $  K_X$ ample
 and $X$ embedded by a suitable multiple of $ K_X$,
 one  has  
$\mathsf P_X\simeq \operp^s U \operp \operp ^t E_8(-1)$ since the index is negative by the index formula (cf. Theorem~\ref{thm:Index}).
The Noether inequality \cite[Theorem VII.3.1]{4authors} stating that $p_g\le \half c_1^2+2$ implies that $p_g\le 2$. 
Furthermore, in case $q>0$, O.\ Debarre \cite{deb1,deb2} has show that $2p_g\le c_1^2$ so that $p_g=0$  in the present situation.  From this and the Noether formula~\eqref{eqn:Noether}, one arrives at the following
sets  of  possible invariants:
\begin{center}
  \begin{tabular}{@{} |c| c|c| c | @{}}
    \hline
    $\chi$  & $1$ & $2$ & $3$ \\ 
    \hline
  $c_2$ & $11$ & $23$  & 35\\ 
\hline
   $(p_g,q)$ & $(0,0) $  & $  (1,0)  $  &$ (2,0)$ \\ 
\hline
$(s,t)$ & $(0,1) $  & $  (2,2) $  &$ (4, 3)$ \\ 
    \hline
  \end{tabular}
\end{center}
 Here are some surfaces within this range of invariants (the list is far from complete!):

\begin{itemize}
\item The  so-called Godeaux-type surfaces, i.e. those with $p_g=q=0$ and $c_1^2=1$. For concrete examples, see e.g.\ \cite[\S VII.10]{4authors}. 
Here $\mathsf P_X$ is unimodular and negative definite of rank $8$. It is known that
then $\mathsf P_X\simeq E_8(-1)$ (cf.\ Table~\ref{tab1}).
\item The Kynev surface from \cite{kynev, todorov2} with $c_1^2=1$, $p_g=1$ and $q=0$ (so that $c_2=23$). 
\item E.\ Horikawa's (simply connected) surface from \cite{hor1} with $c_1^2=1$, $c_2=35$ (so that $b_2(X)=33$).
\end{itemize}

\item The simplest non-unimodular $\mathsf P_X$ are obtained for surfaces $X$ with $c_1^2(X)=2$ and $  K_X$ ample
 and $X$ embedded by a suitable multiple of $ K_X$. Here $\disc{\mathsf P_X}=\pm 2$. As before, using Noether's inequality, Debarre's inequality and the Noether formula,  one arrives at the following
sets  of  possible invariants:

\begin{center}
  \begin{tabular}{@{} |c| c|c| c |c| @{}}
    \hline
    $\chi$  & $1$ & $2$ & $3$& $4$  \\ 
    \hline
  $c_2$ & $10$ & $22$  & $34$ & $46$ \\ 
\hline
   $(p_g,q)$ & $(0,0), (1,1) $  & $  (1,0)  $  &$ (2,0)$ &$(3,0)$\\ 
\hline
  \end{tabular}
\end{center}

These surfaces are known to exist. I give some examples:

\begin{itemize}

\item The numerical Campedelli surfaces, i.e., those with $p_g=q=0$ and $c_1^2=2$.  Again, for examples see e.g.\ \cite[\S VII.10]{4authors} 
For these, $\mathsf P_X$ is negative definite of rank $2$ and with $\disc{\mathsf P_X}=-2$. It is known that $\mathsf P_X\simeq E_7(-1)$. See Remark~\ref{rmk:OnIndefLats} 
and Table~\ref{tab1} below.

\item The surfaces with $p_g=q=1$ and $c_1^2=2$ have been completely classified. See \cite{cat}. Here $\mathsf P_X$ has signature $(2,9)$ and discriminant $-2$.
Such a lattice is isometric to $ \la 2\ra \operp U\operp E_8(-1)$. This follows from Theorem~\ref{thm:kneser} since the given lattice has the
correct signature and discriminant form.
 
 \item Horikawa's surface with $c_1^2=2, p_g=3, q=0$ (and $c_2=46$)  from \cite{hor2}. Here $\mathsf P_X$ has signature $(6,37)$ and discriminant $-2$.
 Such a lattice is isometric to $ \la 2\ra \operp \operp^{ 5} U \operp \operp ^4 E_8(-1)$.
\end{itemize}

\item Let $X$ be an Enriques surface. Then $\mathsf H_X\simeq U\operp E_8(-1)$. Let $c$ be a primitive vector in the $U$ component, say
 $c=  e + f$ where  $\set{e,f}$ is the standard basis of $U$. Then $c^\perp \simeq \la -2 \ra \operp  E_8(-1)$.
 By the main theorem in \cite{wats5} this lattice has class number $1$.
By loc. cit.   for vectors of the form  $c'=d.e+f$, $d\not=\pm 1$,  the class number 
of the lattice $(c')^\perp$ is larger than $1$.

 In fact, to interpret Watson's results, one has to be careful 
since his  terminology differs form what is nowadays usual. 
First of all, Watson only considers quadratic forms and so the associated bilinear forms (the polar forms) 
are  always even. 
His notation  compares to the one used in this note as follows: 
$P=U$, $Q=\la 2\ra$, $B=A_2$. $E=E_8$ so that the two forms of rank $9$
having class number $1$ are $F_9= E_8\operp \la 2 \ra$ and $G_9$, an indecomposable
 form of discriminant $8$ (in loc. cit. the discriminant of forms of odd rank have been divided by $2$). The last form is not isometric to $E_8\operp\la 8\ra$ since $(G_9)_2= \operp^3 U \operp A_2 \operp \la- 3. 2^3\ra$. 
\end{enumerate}

  \end{exmples}

 \begin{rmk} \label{rmk:OnIndefLats}  If $\mathsf P_X$ is definite,  Theorem~\ref{thm:kneser}  does not apply.
However, there are     lists  of low rank definite  lattices that have one isometry class in its genus. See e.g., \cite{wats,wats1,wats2,wats3,wats4,wats5,wats6,wats7}. This leads to the following table.

\begin{table}[htp]
\caption{List of lattices $\mathsf P_X$ for $X$ canonically polarized with  $\chi(X)=12$.} \label{tab1} 
\begin{center}
\begin{tabular}{|c|c|c|}
\hline
$c_1^2,\rank {P_X}$ & lattice & discrim. form\\
\hline
$(1,8)$& $E_8(-1)$ & $0$\\
$(2,7)$& $E_7(-1)$& $\la -1/2\ra$ \\
$(3,6)$& $E_6(-1) $ & $\la 1/3\ra$ \\
$(4,5)$& $D_5(-1) $ & $\la -1/4\ra$\\
$(5,4)$ & $A_4(-1) $ & $\la -4/5\ra$\\
$(6,3$)& $A_2(-1)\operp \la -2\ra$ & $\la 1/3\ra\operp \la -1/2\ra$ \\
$(7,2)$&$\begin{pmatrix} -4&1\\1&-2\end{pmatrix}$ & $\la  1/7\ra $\\
$(8,1)$&$\la -8\ra$  & $\la-1 / 8\ra$\\
\hline
\end{tabular}
\end{center}
\label{tab:pg=0}
\end{table}%
That  the given lattices of rank $8$, $2$ and $1$ have class number $1$ is trivial or else well known. For other ranks I refer to the cited articles by G.\ Watson.  The  lattices in the table indeed have   rank $9-k$ and 
discriminant group $\bZ/ k\bZ$, $k=c_1^2$  and so these  match with those for which  the results  in loc.\ cit.
show that the class number of the  genus equals one.

 \end{rmk}


 \bibliographystyle{acm}

\end{document}